\documentclass[pdflatex,sn-mathphys-num]{sn-jnl}

\usepackage{graphicx}
\usepackage{multirow}
\usepackage{amsmath,amssymb,amsfonts}
\usepackage{amsthm}
\usepackage{mathrsfs}
\usepackage[title]{appendix}
\usepackage{xcolor}
\usepackage{textcomp}
\usepackage{manyfoot}
\usepackage{booktabs}
\usepackage{longtable}
\usepackage{xltabular}
\usepackage[load-configurations=version-1]{siunitx}
\usepackage[algo2e,ruled]{algorithm2e}
\usepackage{pgfplots}
\usepgfplotslibrary{groupplots}
\pgfplotsset{compat=1.18}

\DeclareMathOperator*{\argmin}{argmin}

\newcommand{\plotTitleSixtyFourFull}{$\lvert V\rvert = 64\text{, } \lvert E\rvert=1987$}
\newcommand{\plotTitleSixtyFourMid}{$\lvert V\rvert=64\text{, } \lvert E\rvert=1713$}
\newcommand{\plotTitleSixtyFourSparse}{$\lvert V\rvert=64\text{, } \lvert E\rvert=986$}
\newcommand{\plotTitleOneTwentyEightFull}{$\lvert V\rvert = 128\text{, } \lvert E\rvert=6432$}
\newcommand{\plotTitleOneTwentyEightMid}{$\lvert V\rvert=128\text{, } \lvert E\rvert=4786$}

\begin{document}

\title[HUBO vs QUBO]{Higher-Order vs. Quadratic Binary Optimization: Which Is Better for Probability Optimization with Tensor Sampling?}

\author*[1,2,3]{\fnm{Hadi} \sur{Salloum}}\email{h.salloum@innopolis.ru}
\author[1]{\fnm{Kirill} \sur{Novoselov}}\email{k.novoselov@innopolis.university}
\author[1]{\fnm{Aleksandr} \sur{Pochtarev}}\email{a.pochtarev@innopolis.university}
\author[1,2]{\fnm{Yaroslav} \sur{Kholodov}}\email{ya.kholodov@innopolis.ru}

\affil*[1]{\orgname{Research Center of the Artificial Intelligence Institute, Innopolis University}, \orgaddress{\city{Innopolis}, \country{Russia}}}
\affil[2]{\orgname{Moscow Institute of Physics and Technology}, \orgaddress{\city{Moscow}, \country{Russia}}}
\affil[3]{\orgname{Q Deep}, \orgaddress{\city{Innopolis}, \country{Russia}}}
\abstract{This article explores the comparative strengths of Higher-Order Unconstrained Binary Optimization (HUBO) and Quadratic Unconstrained Binary Optimization (QUBO) models in the context of probability optimization using tensor sampling techniques. HUBO can represent interactions beyond pairwise terms without introducing auxiliary variables, while QUBO benefits from a simpler algebraic structure and a wide range of established solvers. We combine a theoretical analysis of formulation expressiveness and reduction overhead with empirical experiments on synthetic polynomial objectives, RSA factorization, and Max-Cut instances. The results show that native HUBO formulations are often preferable when the original problem contains high-order interactions, because reducing such problems to QUBO increases the effective dimension seen by tensor samplers. At the same time, QUBO remains competitive for naturally quadratic problems. These findings provide practical guidance for selecting an optimization framework that balances model fidelity, dimension growth, and solver performance.}

\maketitle
\section{Introduction}
\label{sec:intro}

Combinatorial optimization focuses on finding the best solution from a finite but often exponentially large set of possibilities. Within this domain, \emph{binary optimization}, where decision variables take values in $\{0,1\}$, plays a central role. Applications span scheduling, network design, facility location, portfolio selection, feature selection, and machine learning \citet{lucas2014ising, jiang2023classifying, marzec2016portfolio, salloum2025performance}. Since these problems are typically NP-hard, research has explored both exact approaches and heuristic algorithms to achieve a balance between computational tractability and solution quality.

Two standard formulations are \textbf{Quadratic Unconstrained Binary Optimization (QUBO)} and its generalization, \textbf{Higher-Order Unconstrained Binary Optimization (HUBO)}. QUBO restricts the objective to linear and quadratic terms, benefiting from an extensive ecosystem of solvers, both classical and quantum \citet{glover, salloum2025mini, salloum2025quantum, salloum2024enhancing, salloum2024quantum}. HUBO, in contrast, allows higher-order interactions among three or more variables, offering richer modeling fidelity without introducing auxiliary variables. However, higher-order formulations may increase complexity and hinder scalability. Thus, the choice between HUBO and QUBO reflects a fundamental trade-off between expressiveness and efficiency.

A promising direction for such problems is \emph{optimization with tensor train sampling}. Representing a probability tensor as a tensor train (TT) enables efficient probabilistic sampling in high-dimensional binary spaces without exhaustive enumeration. Because tensors naturally encode multi-way dependencies, they align well with HUBO's higher-order interactions \citet{TT, dolgov2020approximation}. The PROTES framework \citet{PROTES} demonstrates how tensor-based sampling achieves competitive performance in black-box optimization tasks \cite{salloum2025black, salloumperformance, salloumtensor}. This makes tensor sampling a useful setting in which to ask whether the algebraic simplicity of QUBO compensates for the dimension growth caused by reducing higher-order terms.

These observations lead to a central investigation:

\begin{quote}
\textbf{What are the trade-offs between HUBO and QUBO formulations in probabilistic tensor sampling optimization, and how do these trade-offs evolve when the number of variables is fixed but the order of interactions increases?}
\end{quote}

\subsection{Contributions}
We formalize probabilistic binary optimization under HUBO and QUBO, highlighting structural trade-offs and reformulation costs. We conduct a systematic study of fixed-variable problems with varying interaction orders to isolate higher-order effects. Then we evaluate tensor sampling methods under both native and transformed formulations, comparing efficiency, accuracy, and robustness. Building on these results, we identify when HUBO's expressiveness reduces the effective search dimension and when QUBO's simpler form remains a more natural representation.

\subsection{Paper Structure}
The paper is organized as follows: Section \ref{sec:background} reviews the necessary background on HUBO, QUBO, and tensor-based optimization methods. Section \ref{sec:theory} compares HUBO and QUBO formulations theoretically. Section \ref{sec:experiments} describes the experimental setup and presents the results. Section \ref{sec:discussion} discusses our results, and Section \ref{sec:conclusion} concludes with future directions.

\section{Background}
\label{sec:background}

\subsection{HUBO and QUBO Formulations}

A Higher-Order Unconstrained Binary Optimization (HUBO) objective of order $k$ is
\begin{equation}
  E_{\text{HUBO}}(\mathbf{x}) = \sum_{p=1}^k \sum_{i_1<\dots<i_p} w_{i_1\dots i_p} \prod_{t=1}^p x_{i_t},
\end{equation}
where $w_{i_1\dots i_p}$ are coefficients describing interactions among $p$ binary variables.  
The order $k$ of the HUBO corresponds to the polynomial degree of the objective.

The Quadratic Unconstrained Binary Optimization (QUBO) problem is the special case $k=2$:
\begin{equation}
  E_{\text{QUBO}}(\mathbf{x}) = \sum_i a_i x_i + \sum_{i<j} b_{ij} x_i x_j,
\end{equation}
where $a_{i}$ are linear coefficients and $b_{ij}$ capture pairwise interactions.  
In practice, QUBO constraints are absorbed into the objective via penalty terms, and the objective can be expressed in matrix form:
\begin{equation}
  E_{\text{QUBO}}(\mathbf{x}) = \mathbf{x}^T Q \mathbf{x}.
\end{equation}

Both HUBO and QUBO objectives define discrete polynomial optimization problems over $\mathbf{x} \in \{0,1\}^d$.  
Directly enumerating all $2^d$ configurations is infeasible for large $d$, motivating specialized optimization strategies.

\subsection{Tensor Representations of Probability Tensors}

A binary optimization problem can be associated with a tensor representation.  
For HUBO/QUBO themselves, one could in principle construct an order-$d$ tensor
\[
\mathcal{E} \in \mathbb{R}^{\underbrace{2 \times \dots \times 2}_{d}},
\]
where $\mathcal{E}(i_1,\dots,i_d)$ stores the objective value at $\mathbf{x}=(i_1,\dots,i_d)$.  
However, storing and manipulating this tensor is intractable for large $d$ due to exponential growth.

 Instead, tensor decompositions are applied to the \emph{probability tensor}, which encodes a distribution over binary solutions.

\subsection{Tensor Train Format}
\label{sec:tt_format}

In the TT format, a high-dimensional tensor is represented as a chain of 3D cores:
\begin{equation}
\mathcal{T}[n_1, n_2, \dots, n_d] = 
\sum_{r_0=1}^{R_0} \cdots \sum_{r_d=1}^{R_d} 
\mathcal{G}_1[r_0, n_1, r_1] 
\mathcal{G}_2[r_1, n_2, r_2] \dots 
\mathcal{G}_d[r_{d-1}, n_d, r_d],
\end{equation}
where $\mathcal{G}_k \in \mathbb{R}^{R_{k-1} \times N_k \times R_k}$ are TT-cores, and $R_k$ are TT-ranks.  
This representation reduces memory requirements from $O(2^d)$ to 
$O(d \cdot N \cdot R^2)$, where $N = \max_k N_k$ and $R = \max_k R_k$.  
Basic linear algebra operations (norms, contractions, etc.) can be performed in polynomial time in $d$.

\subsection{Tensor Train in Optimization}
\label{sec:tensor_sampling}

The key idea is that the probability distribution over solutions $\mathcal{P}$ can be stored in TT-format.  
This makes it possible to efficiently sample solutions and update the distribution in high dimensions.  
Two leading tensor-based optimization methods are:

\begin{enumerate}
    \item \textbf{PROTES} \citet{PROTES}: A probabilistic method that parameterizes $\mathcal{P}$ in TT-format, samples candidate solutions, evaluates them, and updates $\mathcal{P}$ to bias toward high-quality solutions. See Appendix~\ref{app:complexity} for complexity analysis.
    \item \textbf{TTOpt} \citet{TTOpt}: A deterministic method based on the maximum-volume principle. It identifies representative submatrices during tensor unfolding and uses them to approximate the maximum entry of $\mathcal{P}$. See Appendix~\ref{app:complexity}.
\end{enumerate}

These methods are not applied to the HUBO/QUBO tensor directly, but rather to the \emph{probability tensor} that guides the search over the exponentially large solution space.  

A practical comparison of scalability is shown in Table~\ref{tab:tensor_limits}.

\begin{table}[h!]
\centering
\caption{Scalability of tensor-based optimization methods (see Appendices \ref{app:tt_optimization} and \ref{app:complexity}). \citet{PROTES, TTOpt, OptimaTT}}
\label{tab:tensor_limits}
\begin{tabular}{l c}
\hline
\textbf{Method} & \textbf{Maximum dimension $d$} \\
\hline
PROTES & $\sim 1000$ \\
TTOpt & $\sim 500$ \\
Optima-TT & $\sim 100$ \\
\hline
\end{tabular}
\end{table}

\section{Theoretical Comparison}
\label{sec:theory}
In different combinatorial optimization problems, the choice between HUBO and QUBO can be more natural depending on the structure of the objective. In this section, we compare the two formulations from the point of view of expressiveness and reduction overhead. The main issue is not only whether a high-order polynomial can be quadratized, but also how many additional variables are introduced and how this changes the tensor dimension handled by the sampler.

\subsection{Expressiveness}

In the context of combinatorial optimization, \emph{expressiveness} refers to the ability of a problem formulation to accurately represent interactions between variables. QUBO models are limited to pairwise interactions, i.e., terms of the form $x_i x_j$, whereas HUBO formulations can include higher-order terms $x_{i_1} x_{i_2} \dots x_{i_k}$ for $k > 2$ and therefore represent dependencies among several objects directly. We provide examples to illustrate when a QUBO formulation arises more naturally, and when a HUBO formulation is more appropriate.

\paragraph{Max-Cut Problem.}

Given an undirected graph $G(V, E)$ with a set of vertices $V$ and a set of edges $E$, the Max-Cut problem asks for a partition of $V$ into two sets such that the number of edges crossing the partition is maximized.

We can formulate this problem as a QUBO by introducing binary variables such that $x_j=1$ if vertex $j$ is in the first set and $x_j=0$ otherwise. The expression $x_i+x_j-2x_ix_j$ is equal to 1 only if exactly one of the two variables is equal to 1, and it is equal to 0 in all other cases \citet{glover}. Thus, this expression indicates whether the endpoints of an edge belong to different parts of the cut.

So, the Max-Cut problem can be formulated as: 

\begin{equation} \label{eq:max-cut first}
     y = \max \sum_{(i,j) \in E} \left( x_i + x_j - 2x_i x_j \right)
\end{equation}

We see that this is a quadratic form, so we can reformulate this problem as a QUBO: 

\begin{equation}
    y = \max_{\mathbf{x} \in \{0, 1\}^n} \mathbf{x}^TQ\mathbf{x}
\end{equation}

The condition of a graph’s vertex set, such as the existence of an edge, is naturally formulated for two vertices; therefore, the QUBO formulation is the most interpretable in this case.

\paragraph{RSA Problem.}
Let a large composite number $N$ be given.
The RSA problem consists of representing this number as the product of two prime numbers $$N = p \times q$$

We can easily reduce this to an optimization problem \citet{rsa}:

\begin{equation}
\argmin_{p,q \in \mathbf{N}} \|pq - N\|_2 \label{eq:rsa}
\end{equation}

\begin{equation}
\|pq - N\|^2_2 = p^2q^2 - 2pqN + N^2
\end{equation}

Next, we will quantize the numbers $p$ and $q$ and substitute them into the objective functional:

\begin{equation}
p = \sum_{l=0}^{n-1} 2^l q_l, \quad q = \sum_{l=0}^{n-1} 2^l q_{n+l}.
\end{equation}

\begin{align}
p^2 q^2 &= \left( \sum_{l=0}^{n-1} 2^l q_l \right)^2 \left( \sum_{l=0}^{n-1} 2^l q_{n+l} \right)^2 \nonumber \\
&= \left( \sum_{l=0}^{n-1} 2^{2l} q_l + \sum_{l_1 < l_2} 2^{l_1 + l_2 + 1} q_{l_1} q_{l_2} \right) \left( \sum_{l=0}^{n-1} 2^{2l} q_{n+l} + \sum_{l_1 < l_2} 2^{l_1 + l_2 + 1} q_{n+l_1} q_{n+l_2} \right) \nonumber \\
&= \sum_{l_1=0}^{n-1} \sum_{l_2=0}^{n-1} 2^{2(l_1 + l_2)} q_{l_1} q_{n+l_2} + \sum_{l_1 < l_2} \sum_{l_3 < l_4} 2^{2l_1 + l_2 + l_3 + 1} (q_{l_1} q_{n+l_2} q_{n+l_3} + q_{l_2} q_{l_3} q_{n+l_4}) \nonumber \\
&\quad + \sum_{l_1 < l_2} \sum_{l_3 < l_4} 2^{l_1 + l_2 + l_3 + l_4 + 2} q_{l_1} q_{l_2} q_{n+l_3} q_{n+l_4}.
\end{align}

\begin{equation}
-2pqN = -2N \left( \sum_{l=0}^{n-1} 2^l q_l \right) \left( \sum_{l=0}^{n-1} 2^l q_{n+l} \right) = \sum_{l_1=0}^{n-1} \sum_{l_2=0}^{n-1} (-2^{l_1 + l_2 + 1} N q_{l_1} q_{n+l_2}).
\end{equation}

Minimising the function in this form provides an equivalent reformulation of the RSA problem as a HUBO problem.

In addition to their natural interpretability, HUBO and QUBO formulations often differ in the order of the tensor to which they can be reduced.
Recall that HUBO and QUBO can be expressed as polynomial objectives in $d$ binary variables, with degree two for QUBO. As discussed in Section~\ref{sec:background}, the probability tensor used by tensor sampling methods has one mode per binary variable. Therefore, a reduction that increases the number of variables directly increases the tensor dimension, even if the polynomial degree becomes lower. This makes it critical to determine how many variables $x_i$ are required by each formulation and what practical tensor-size limits are imposed by the chosen sampling method (Table~\ref{tab:tensor_limits}).

\subsection{Reduction of HUBO to QUBO}
A typical reduction introduces auxiliary variables $y$ and penalty terms to enforce $y = \prod x_{i_t}$. The overhead grows with the number and order of monomials, and this overhead is particularly important for tensor-based methods because each auxiliary variable adds a new tensor mode.

As was shown in \citet{rsa}, we can reduce the order-3 monomial $cxyz$, where $c \in \mathbb{R}$ in the following way:

\begin{equation} \label{eq:3to2}
c x y z = 
\begin{cases} 
c w(x + y + z - 2), & c < 0 \\
c \left\{ w(x + y + z - 1) + (xy + yz + zx) - (x + y + z) + 1 \right\}, & c > 0 
\end{cases}
\end{equation}

Continuing this idea, we can consider higher-degree monomials by temporarily grouping variables into a single one to obtain a cubic term.
For example, we can reduce the monomial $a_1 a_2 b_1 b_2$ with a positive coefficient as follows:
\begin{equation}
\begin{aligned}
a_1 a_2 b_1 b_2 &= a_1 \left( x_1 a_2 + x_1 b_1 + a_2 b_1 + a_2 b_2 + b_1 b_2 - a_2 - b_2 - x_1 + 1 \right) \\
&= x_2 x_1 + x_2 a_1 + x_1 a_2 + x_1 a_1 + x_1 a_2 - x_1 - a_1 - a_2 - x_2 + 1 \\
&\quad + x_3 x_1 + x_3 a_1 + x_1 b_1 + x_1 b_1 + a_1 b_1 - x_1 - a_1 - b_1 - x_3 + 1 \\
&\quad + x_4 x_1 + x_4 a_1 + x_1 b_2 + x_1 b_2 + a_1 b_2 - x_1 - a_1 - b_2 - x_4 + 1 \\
&\quad + x_5 x_1 + x_5 a_1 + a_1 a_2 + a_2 b_1 - a_1 - a_2 - b_1 - x_5 + 1 \\
&\quad + x_6 x_1 + x_6 a_1 + a_1 a_2 + a_2 b_2 - a_1 - a_2 - b_2 - x_6 + 1 \\
&\quad + x_7 x_1 + x_7 a_1 + a_1 b_1 + a_1 b_2 + b_1 b_2 - a_1 - b_1 - b_2 - x_7 + 1 \\
&\quad - a_1 a_2 - a_1 b_1 - a_1 b_2 - a_1 x_1 + a_1
\end{aligned}
\end{equation}
We observe that reducing a cubic monomial requires the introduction of one auxiliary variable, whereas a quartic monomial may require as many as seven. As the degree of the monomial increases, the number of auxiliary variables introduced by this reduction method grows quickly. Given the limitations of tensor sampling methods, this is critical for our purposes: a HUBO-to-QUBO reduction may produce a formally quadratic model that is nevertheless harder for tensor sampling because of the resulting increase in dimension.

\subsection{Max-Cut formulation as a HUBO problem}
Recall that, as a consequence of expression \ref{eq:3to2}, for $x_{0}, x_{1}, x_{2}, x_{3} \in \{0,1\}$ and $c \in \mathbb{R}$ with $c<0$, the following equality holds:
\begin{equation} \label{eq:max-cut1}
c x_{1} x_{2} x_{3} \;=\; c x_{0}x_{1} + c x_{0}x_{2} + c x_{0}x_{3} - 2c x_{0}
\end{equation}

Thus, when $c=-2$, the following equality holds:
\begin{equation} \label{eq:max-cut2}
    -2 x_{1} x_{2} x_{3} - 4x_{0} \;=\; -2 x_{0}x_{1} -2 x_{0}x_{2}  -2 x_{0}x_{3}
\end{equation}

This means that if we can identify terms of the form $-2x_{0}x_{1} - 2x_{0}x_{2} - 2x_{0}x_{3}$ in our QUBO formulation \ref{eq:max-cut first}, 
then we can equivalently replace them with cubic and linear terms by increasing the order of the problem.

To this end, we modify the formulation of the Max-Cut problem. Let our graph $G = (V, E)$ be defined on $n$ vertices. 
Then $E_{\text{full}}$ denotes the set of all possible edges $(i, j)$ with $0 \leq i < j < n$ in this graph. 
Let $M = E_{\text{full}} \setminus E$ be the set of edges that are missing from our graph. 

First, we compute the objective function \ref{eq:max-cut first} for the case of a complete graph. Then, we successively take quadruples of coordinates of the vector $x$, namely $x_{i}, x_{i+1}, x_{i+2}, x_{i+3}$, 
and replace their occurrence in the objective function in the form 
$-2x_{i}x_{i+1} - 2x_{i}x_{i+2} - 2x_{i}x_{i+3}$
with cubic and linear terms using expression \ref{eq:max-cut2}. 

After that, we subtract from the obtained result the value of the objective function computed on the set $M$. The resulting polynomial is equivalent to the original Max-Cut objective, but it contains cubic terms introduced by the replacement step.
We illustrate this procedure in Algorithm~\ref{alg:Max-Cut_hubo}.

\begin{algorithm2e}
\caption{Sparse Graph Max-Cut as 3rd-order HUBO}
\label{alg:Max-Cut_hubo}

\KwIn{Graph $G = (V, E)$ and Cubic replacement formula: $-2 x_{1} x_{2} x_{3} - 4x_{0} = -2 x_{0}x_{1} -2 x_{0}x_{2}  -2 x_{0}x_{3}$}
\KwOut{Cubic polynomial $Q_{sparse,HUBO}$ representing Max-Cut}

$E_{full} \leftarrow$ all pairs of nodes in $V$\;
$M \leftarrow E_{full} \setminus E$ \tcp*{Missing edges}
$Q_{full} \leftarrow 0$\;

\For{each edge $(i,j) \in E_{full}$}{
  $Q_{full} \leftarrow Q_{full} + (x_i + x_j - 2x_i x_j)$\;
}

\For{$i \leftarrow 0$ \KwTo $n$ \textbf{step} $4$}{
  $Q_{full} \leftarrow Q_{full} + 2x_ix_{i+1}+ 2x_ix_{i+2}+ 2x_ix_{i+3}$\;
  $Q_{full} \leftarrow Q_{full} - 2 (x_{i+1} x_{i+2} x_{i+3} + 2 x_i)$\;
}

$Q_{sparse,HUBO} \leftarrow Q_{full} - \sum_{(i,j)\in M}(x_i + x_j - 2x_i x_j)$\;

\KwRet $Q_{sparse,HUBO}$\;

\end{algorithm2e}

\textbf{Theorem 1:}
Algorithm~\ref{alg:Max-Cut_hubo} is correct.

\textit{Proof}
Indeed, since we only performed equivalent substitutions, even after increasing the order of the problem, we can write $Q_{\text{full}}$ in the form
\[
Q_{\text{full}} = \sum_{(i,j) \in E_{\text{full}}} \bigl(x_{i} + x_{j} - 2x_{i}x_{j}\bigr)
= \sum_{(i,j) \in M} \bigl(x_{i} + x_{j} - 2x_{i}x_{j}\bigr)
+ \sum_{(i,j) \in E} \bigl(x_{i} + x_{j} - 2x_{i}x_{j}\bigr).
\]

Therefore, in the final step of the algorithm we have
\[
Q_{\text{sparse}}
= \sum_{(i,j) \in M} \bigl(x_{i} + x_{j} - 2x_{i}x_{j}\bigr)
+ \sum_{(i,j) \in E} \bigl(x_{i} + x_{j} - 2x_{i}x_{j}\bigr)
- \sum_{(i,j) \in M} \bigl(x_{i} + x_{j} - 2x_{i}x_{j}\bigr).
\]

This simplifies to
\[
Q_{\text{sparse}} = \sum_{(i,j) \in E} \bigl(x_{i} + x_{j} - 2x_{i}x_{j}\bigr),
\]
which is exactly the original Max-Cut objective. 
$\square$

\section{Experimental Setup and Results}
\label{sec:experiments}
We conducted experiments to evaluate tensor sampling performance on HUBO, QUBO, and transformed HUBO-to-QUBO or QUBO-to-HUBO formulations. Test cases included arbitrary polynomial functions, RSA factorization problems, and Max-Cut instances on random graphs. For each problem, we measured solution quality and computation time or evaluation budget across different problem sizes ($n$) and polynomial orders ($k$).

Both PROTES and TTOpt solvers were applied where appropriate to examine the trade-off between solution quality and computational budget, and to verify whether the observed results generalize across different optimization methods. For transformed QUBO formulations, auxiliary variables introduced by the reduction were included in the search space; therefore, the comparison reflects the effective dimension seen by the solver rather than only the polynomial degree of the objective. Convergence and best-found solutions were tracked over iterations, and the Max-Cut results are presented in Tables~\ref{tab:maxcut_protes} and \ref{tab:maxcut_ttopt}, as well as in Figures~\ref{fig:experiments64} and~\ref{fig:experiments128}, for comparison.\footnote{We removed Optima-TT because it can only solve small problems; see Table~\ref{tab:tensor_limits}.}

\subsection{Arbitrary Functions}
We first compare HUBO and HUBO-to-QUBO formulations on synthetic polynomial objectives. These functions are small enough to make the effect of the reduction transparent, but they already contain cubic, quartic, and fifth-order interactions:
\begin{equation}
    f_1(\mathbf{x})= -3x_1x_3x_6 + 6x_2x_1x_3 - x_3x_4x_5x_6 + x_3x_6 + 2x_2x_5 - x_5 - x_1
\end{equation}
\begin{equation}
    f_2(\mathbf{x})= 5x_1x_2x_3 - 4x_2x_3x_4 - 2x_1x_3x_4 + x_1x_3 + x_2 - x_4
\end{equation}
\begin{equation}
    f_3(\mathbf{x})= -2x_1x_2x_3x_4x_5 - x_1x_3x_4x_5x_6 + 2x_1x_2x_3x_6 + x_2x_3x_4x_5 - 3x_1x_5 + x_2x_5
\end{equation}
\begin{multline}
f_4(\mathbf{x}) = -2x_1x_2x_3x_4x_5 + x_2x_3x_5x_6x_7 - 3x_2x_3x_4x_5 + x_1x_2x_3x_4 \\
 - x_3x_5x_6x_7 + 2x_1x_2x_3 - 2x_2x_3x_6 + x_1x_6 - x_6x_7 - x_7
\end{multline}

\begin{table}[h]
  \centering
  \caption{HUBO and HUBO-to-QUBO quality and time comparison}
  \begin{tabular}{|c|c|c|c|c|c|c|}
    \hline
    Function & $n$ & order $k$ & HUBO result & HUBO-to-QUBO result & HUBO time & HUBO-to-QUBO time \\ \hline
    $f_1(\mathbf{x})$ & 6 & 4 & -5 & -5 & 1.8s & \textbf{1.1s} \\ \hline
    $f_2(\mathbf{x})$ & 4 & 3 & -4 & -4 & 1.7s & \textbf{1.1s} \\ \hline
    $f_3(\mathbf{x})$ & 6 & 5 & -4 & -4 & 1.8s & \textbf{1.2s} \\ \hline
    $f_4(\mathbf{x})$ & 7 & 5 & -7 & -7 & 1.8s & \textbf{1.5s} \\ \hline
  \end{tabular}
  \label{tab:blowup}
\end{table}
For the first four functions, both formulations return the same objective value, although the reduced QUBO is sometimes faster because the problem size is still small. However, this behavior does not persist for all high-order objectives. For example, for the third-degree function
\begin{equation*}
f_5(\mathbf{x}) =
    \sum_{i=0}^{n-1} \sum_{j=i+1}^{n-1} \sum_{\substack{k=j+1 \\ \text{$k$ is odd}}}^{n-1} 
(-1)^j (j - i)  x_i  x_j  x_k
\end{equation*}
the output and time become different for some values of $n$.
\begin{table}[h]
  \centering
  \caption{HUBO and HUBO-to-QUBO for $f_5(\mathbf{x})$}
  \begin{tabular}{|c|c|c|c|c|c|}
    \hline
    $n$ & order $k$ & HUBO result & HUBO-to-QUBO result & HUBO time & HUBO-to-QUBO time \\ \hline
    10 & 3 & \textbf{-62} & -60 & 7.7s & 2.5s \\ \hline
    16 & 3 & \textbf{-462} & -392 & 33s & 17s \\ \hline
    18 & 3 & \textbf{-761} & -682 & 58s & 75s \\ \hline
  \end{tabular}
  \label{tab:3}
\end{table}

For fourth-degree functions, the situation is more complex because degree reduction with a positive coefficient adds seven auxiliary variables in the construction considered above, and a negative coefficient also increases the number of variables. Consider the function
\begin{equation*}
f_6(\mathbf{x})=
    \sum_{i=0}^{n-1} 
\sum_{j=i+1}^{n-1} 
\sum_{\substack{k=j+1 \\ k\equiv j+1 \pmod{2}}}^{n-1} 
\sum_{\substack{p=k+1 \\ p\equiv k+1 \pmod{2}}}^{n-1} 
(-1)^p (j - i + k)  x_i  x_j  x_k  x_p
\end{equation*}
After the HUBO-to-QUBO transformation, PROTES does not recover the correct answer for the tested cases, while the native HUBO formulation does. This indicates that the additional variables introduced by the reduction can dominate the apparent advantage of a quadratic objective.
\begin{table}[h]
  \centering
  \caption{HUBO and HUBO-to-QUBO for $f_6(\mathbf{x})$}
  \begin{tabular}{|c|c|c|c|c|c|}
    \hline
    $n$ & order $k$ & HUBO result & HUBO-to-QUBO result & HUBO time & HUBO-to-QUBO time \\ \hline
    10 & 4 & \textbf{-193} & -106 & 51s & 20s \\ \hline
    16 & 4 & \textbf{-1957} & 436 & 600s & 4928s \\ \hline
  \end{tabular}
  \label{tab:4}
\end{table}

\subsection{RSA problem}
As shown in \eqref{eq:rsa}, the RSA factorization problem can be reduced to a HUBO optimization problem:
\begin{equation}
\argmin_{p,q \in \mathbf{N}} \|pq - N\|_2 
\end{equation}

\begin{figure}[htbp]
    \centering
    \begin{tikzpicture}
        \begin{groupplot}[
            group style={
                group size=3 by 1,
                horizontal sep=1.5cm,
                vertical sep=1.5cm,
            },
            width=0.3\textwidth,
            height=0.3\textwidth,
            grid=major,
            xlabel={Budget},
            ylabel={$y_{opt}$},
            legend pos=north east,
            legend style={font=\scriptsize},
            every axis plot/.append style={mark=*, mark size=1.5pt, line width=1pt},
            xmin=0, xmax=50000,
            xtick={0,10000,20000,30000,40000,50000},
            xticklabels={$0$, $1e4$,$2e4$,$3e4$,$4e4$,$5e4$},
            xticklabel style={/pgf/number format/fixed},
            scaled x ticks=false,
            tick label style={font=\tiny},
            title style={font=\small},
        ]

        \nextgroupplot[title=\plotTitleSixtyFourFull, ymin=-1026, ymax=-1010]
        \addplot+[blue] coordinates {
            (1000.0, -1016.0) (3000.0, -1017.0) (6000.0, -1017.0) (8000.0, -1017.0) (9000.0, -1019.0) (10000.0, -1020.0)
            (11000.0, -1020.0) (14000.0, -1020.0) (15000.0, -1020.0)
            (16000.0, -1021.0) (18000.0, -1022.0) (21000.0, -1023.0)
            (22000.0, -1023.0) (27000.0, -1023.0) 
            (40000.0, -1023.0) (50000.0, -1023.0)
        };
        \addlegendentry{QUBO}
        \addplot+[orange] coordinates {
            (1000.0, -1012.0) (3000.0, -1012.0) (6000.0, -1013.0) (8000.0, -1014.0) (9000.0, -1017.0) (10000.0, -1017.0) 
            (11000.0, -1019.0) (14000.0, -1020.0) (15000.0, -1021.0)
            (16000.0, -1021.0) (18000.0, -1021.0) (21000.0, -1021.0)
            (22000.0, -1022.0)
            (27000.0, -1023.0) (40000.0, -1024.0) (50000.0, -1024.0)
        };
        \addlegendentry{HUBO}

        \nextgroupplot[title=\plotTitleSixtyFourMid, ymin=-930, ymax=-880]
        \addplot+[blue] coordinates {
            (1000.0, -891.0) (2000.0, -891.0) (3000.0, -897.0) (4000.0, -898.0) (6000.0, -898.0) (7000.0, -903.0) (8000.0, -904.0) (10000.0, -904.0)
            (11000.0, -910.0) (12000.0, -911.0) (13000.0, -911.0) (14000.0, -915.0) (15000.0, -915.0) (16000.0, -917.0)
            (17000.0, -918.0) (19000.0, -919.0) (20000.0, -921.0) (23000.0, -921.0) (24000.0, -922.0) (25000.0, -922.0)
            (28000.0, -923.0) (29000.0, -924.0) (32000.0, -924.0) (50000.0, -924.0)
        };
        \addlegendentry{QUBO}
        \addplot+[orange] coordinates {
            (1000.0, -881.0) (2000.0, -883.0) (3000.0, -887.0) (4000.0, -887.0) (6000.0, -890.0) (7000.0, -896.0) (8000.0, -896.0) 
            (10000.0, -904.0) (11000.0, -904.0) (12000.0, -904.0) (13000.0, -910.0) (14000.0, -910.0) (15000.0, -913.0) (16000.0, -918.0)
            (17000.0, -919.0) (19000.0, -919.0) (20000.0, -924.0) (23000.0, -926.0) (24000.0, -926.0) (25000.0, -927.0)
            (28000.0, -927.0) (29000.0, -927.0) (32000.0, -928.0) (50000.0, -928.0)
        };
        \addlegendentry{HUBO}

        \nextgroupplot[title=\plotTitleSixtyFourSparse, ymin=-575, ymax=-520]
        \addplot+[blue] coordinates {
            (1000.0, -533.0) (2000.0, -533.0) (3000.0, -535.0) (5000.0, -535.0) (7000.0, -545.0) (8000.0, -546.0) (10000.0, -546.0)
            (11000.0, -551.0) (14000.0, -551.0) (15000.0, -551.0) (16000.0, -551.0) (17000.0, -556.0) (19000.0, -556.0) (20000.0, -557.0) (21000.0, -560.0)
            (22000.0, -563.0) (23000.0, -563.0) (24000.0, -563.0) (25000.0, -563.0) (27000.0, -564.0) (28000.0, -566.0) (30000.0, -566.0) (31000.0, -567.0) (34000.0, -567.0)
            (35000.0, -569.0) (45000.0, -571.0) (46000.0, -573.0) (50000.0, -573.0)
        };
        \addlegendentry{QUBO}
        \addplot+[orange] coordinates {
            (1000.0, -524.0) (2000.0, -525.0) (3000.0, -525.0) (5000.0, -539.0) (7000.0, -539.0) (8000.0, -539.0) (10000.0, -542.0) (11000.0, -542.0)
            (14000.0, -543.0) (15000.0, -548.0) (16000.0, -551.0) (17000.0, -554.0)
            (19000.0, -560.0) (20000.0, -560.0) (21000.0, -560.0) (22000.0, -560.0) (23000.0, -563.0) (24000.0, -566.0) (25000.0, -569.0) (27000.0, -569.0)
            (28000.0, -570.0) (30000.0, -571.0) (31000.0, -571.0) (34000.0, -572.0) (35000.0, -574.0) (45000.0, -574.0) (46000.0, -574.0)
            (50000.0, -574.0)
        };
        \addlegendentry{HUBO}

        \end{groupplot}
    \end{tikzpicture}
    \caption{Experiments with $n=64$ Max-Cut problem comparing QUBO and HUBO performance via PROTES.}
    \label{fig:experiments64}
\end{figure}
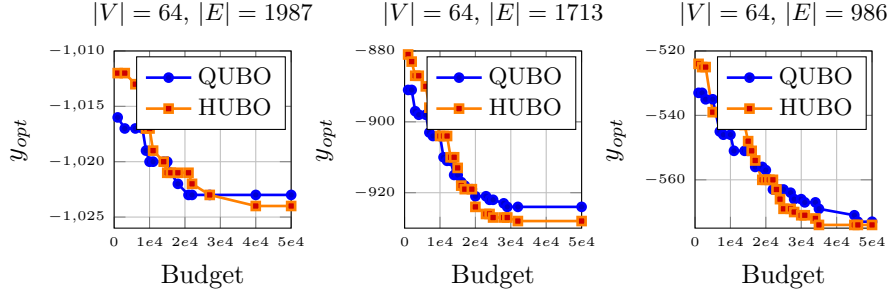

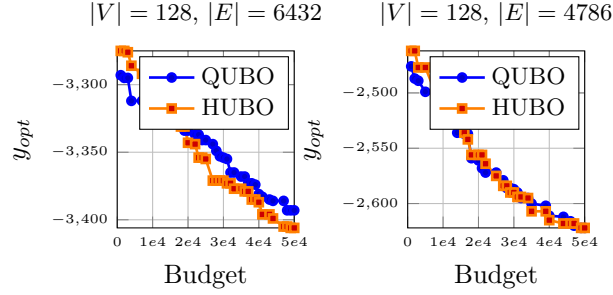
\begin{figure}[htbp]
    \centering
    \begin{tikzpicture}
        \begin{groupplot}[
            group style={
                group size=2 by 1,
                horizontal sep=1.5cm,
                vertical sep=1.5cm,
            },
            width=0.3\textwidth,
            height=0.3\textwidth,
            grid=major,
            xlabel={Budget},
            ylabel={$y_{opt}$},
            legend pos=north east,
            legend style={font=\scriptsize},
            every axis plot/.append style={mark=*, mark size=1.5pt, line width=1pt},
            xmin=0, xmax=50000,
            xtick={0,10000,20000,30000,40000,50000},
            xticklabels={$0$, $1e4$,$2e4$,$3e4$,$4e4$,$5e4$},
            xticklabel style={/pgf/number format/fixed},
            scaled x ticks=false,
            tick label style={font=\tiny},
            title style={font=\small},
        ]

        \nextgroupplot[title=\plotTitleOneTwentyEightFull, ymin=-3406, ymax=-3275]
        \addplot+[blue] coordinates {
            (1000.0, -3293.0) (2000.0, -3295.0) (3000.0, -3295.0) (4000.0, -3312.0) (7000.0, -3312.0) (8000.0, -3314.0) (10000.0, -3314.0) (11000.0, -3314.0) (12000.0, -3314.0) (13000.0, -3314.0) (16000.0, -3321.0) (17000.0, -3328.0) (18000.0, -3331.0) (19000.0, -3334.0) (20000.0, -3334.0) (22000.0, -3337.0) (23000.0, -3337.0) (24000.0, -3341.0) (25000.0, -3341.0) (27000.0, -3344.0) (28000.0, -3349.0) (29000.0, -3353.0) (30000.0, -3354.0) (31000.0, -3355.0) (32000.0, -3365.0) (33000.0, -3365.0) (35000.0, -3368.0) (36000.0, -3368.0) (37000.0, -3373.0) (38000.0, -3373.0) (39000.0, -3374.0) (40000.0, -3381.0) (41000.0, -3383.0) (43000.0, -3385.0) (44000.0, -3386.0) (47000.0, -3386.0) (48000.0, -3393.0) (49000.0, -3393.0) (50000.0, -3393.0)
        };
        \addlegendentry{QUBO}
        \addplot+[orange] coordinates {
            (1000.0, -3275.0) (2000.0, -3275.0) (3000.0, -3276.0) (4000.0, -3286.0) (7000.0, -3292.0) (8000.0, -3292.0) (10000.0, -3297.0) (11000.0, -3306.0) (12000.0, -3313.0) (13000.0, -3314.0) (16000.0, -3327.0) (17000.0, -3327.0) (18000.0, -3331.0) (19000.0, -3331.0) (20000.0, -3343.0) (22000.0, -3344.0) (23000.0, -3354.0) (24000.0, -3354.0) (25000.0, -3355.0) (27000.0, -3371.0) (28000.0, -3371.0) (29000.0, -3371.0) (30000.0, -3371.0) (31000.0, -3373.0) (32000.0, -3373.0) (33000.0, -3377.0) (35000.0, -3377.0) (36000.0, -3379.0) (37000.0, -3379.0) (38000.0, -3385.0) (39000.0, -3385.0) (40000.0, -3387.0) (41000.0, -3396.0) (43000.0, -3396.0) (44000.0, -3399.0) (47000.0, -3405.0) (48000.0, -3405.0) (49000.0, -3406.0) (50000.0, -3406.0)
        };
        \addlegendentry{HUBO}

        \nextgroupplot[title=\plotTitleOneTwentyEightMid, ymin=-2622, ymax=-2462]
        \addplot+[blue] coordinates {
            (1000.0, -2476.0) (2000.0, -2487.0) (3000.0, -2489.0) (5000.0, -2499.0) (8000.0, -2501.0) (10000.0, -2509.0) (11000.0, -2514.0) (12000.0, -2515.0) (14000.0, -2536.0) (16000.0, -2537.0) (17000.0, -2537.0) (18000.0, -2559.0) (20000.0, -2561.0) (21000.0, -2568.0) (22000.0, -2572.0) (25000.0, -2572.0) (27000.0, -2579.0) (28000.0, -2586.0) (29000.0, -2586.0) (30000.0, -2587.0) (31000.0, -2590.0) (32000.0, -2595.0) (34000.0, -2595.0) (35000.0, -2600.0) (39000.0, -2602.0) (40000.0, -2611.0) (44000.0, -2612.0) (46000.0, -2616.0) (47000.0, -2620.0) (49000.0, -2622.0) (50000.0, -2622.0)
        };
        \addlegendentry{QUBO}
        \addplot+[orange] coordinates {
            (1000.0, -2462.0) (2000.0, -2462.0) (3000.0, -2477.0) (5000.0, -2477.0) (8000.0, -2501.0) (10000.0, -2501.0) (11000.0, -2521.0) (12000.0, -2527.0) (14000.0, -2527.0) (16000.0, -2536.0) (17000.0, -2542.0) (18000.0, -2556.0) (20000.0, -2556.0) (21000.0, -2556.0) (22000.0, -2564.0) (25000.0, -2575.0) (27000.0, -2584.0) (28000.0, -2584.0) (29000.0, -2590.0) (30000.0, -2590.0) (31000.0, -2594.0) (32000.0, -2594.0) (34000.0, -2595.0) (35000.0, -2607.0)
            (39000.0, -2607.0) (40000.0, -2615.0) (44000.0, -2618.0) (46000.0, -2618.0) (47000.0, -2618.0) (49000.0, -2622.0) (50000.0, -2622.0)
        };
        \addlegendentry{HUBO}

        \end{groupplot}
    \end{tikzpicture}
    \caption{Experiments with $n=128$ Max-Cut problem comparing QUBO and HUBO performance via PROTES.}
    \label{fig:experiments128}
\end{figure}

For this problem, an objective value of zero corresponds to exact factorization. After reducing the HUBO formulation to QUBO, PROTES does not reach the global minimum in the tested cases, while the native HUBO formulation does. Results for different values of $n$ are shown in Table~\ref{tab:5}.

\begin{table}[t]
  \centering
  \caption{HUBO and HUBO-to-QUBO for RSA problem}
  \begin{tabular}{|c|c|c|c|c|c|}
    \hline
    $n$ & order $k$ & HUBO result & HUBO-to-QUBO result & HUBO time & HUBO-to-QUBO time \\ \hline
    15 & 4 & \textbf{0} & 1 & 6s & 2.9s \\ \hline
    21 & 4 & \textbf{0} & 9 & 4.9s & 5.2s \\ \hline
    35 & 4 & \textbf{0} & 17 & 5s & 18.5s \\ \hline
    77 & 4 & \textbf{0} & 1681 & 14s & 19.5s \\ \hline
  \end{tabular}
  \label{tab:5}
\end{table}

As shown, even for these small examples PROTES fails to obtain the exact answer after the QUBO reduction. In addition, the HUBO formulation is faster for the larger tested instances, which is consistent with the dimension-growth argument.
\subsection{Max-Cut problem}
In Tables~\ref{tab:maxcut_protes} and \ref{tab:maxcut_ttopt} we compare the results of QUBO and HUBO (obtained using Algorithm~\ref{alg:Max-Cut_hubo}) for the Max-Cut problem. The optimization was carried out using PROTES and TTOpt. In addition to the solution found by each method, we also report the budget required to obtain it. For all experiments, random graphs were used, and their characteristics are provided in the tables. Figures~\ref{fig:experiments64} and \ref{fig:experiments128} illustrate the comparison of QUBO and HUBO Max-Cut for $n=64$ and $n=128$.

\begin{table}[t!]
  \centering
  \caption{QUBO and QUBO-to-HUBO for Max-Cut problem via PROTES}
  \label{tab:maxcut_protes}
  \begin{tabular}{|c|c|c|c|c|c|}
    \hline
    $n$ & Number of edges & HUBO result & QUBO result & HUBO budget & QUBO budget \\ \hline
    16 & 48 & 35 & \textbf{36} & 3.0e+03 & \textbf{2.0e+03} \\ \hline
    20 & 76 & \textbf{54} & \textbf{54} & \textbf{6.0e+03} & 1.1e+04 \\ \hline
    20 & 30 & \textbf{25} & \textbf{25} & 8.0e+03 & \textbf{4.0e+03} \\ \hline
    28 & 360 & \textbf{195} & \textbf{195} & \textbf{1.1e+04} & 1.2e+04 \\ \hline
    28 & 121 & \textbf{82} & \textbf{82} & \textbf{2.0e+04} & 3.6e+04 \\ \hline
    28 & 171 & \textbf{110} & \textbf{110} & \textbf{9.0e+03} & 1.1e+04 \\ \hline
    28 & 300 & \textbf{173} & \textbf{173} & 1.0e+04 & \textbf{9.0e+03} \\ \hline
  \end{tabular}
\end{table} 

\begin{table}[t!]
  \centering
  \caption{QUBO and QUBO-to-HUBO for Max-Cut problem via TTOpt}
  \label{tab:maxcut_ttopt}
  \begin{tabular}{|c|c|c|c|c|c|}
    \hline
    $n$ & Number of edges & HUBO result & QUBO result & HUBO evals & QUBO evals \\ \hline
    16 & 48 & \textbf{35} & \textbf{35} & \textbf{1.51e+04} & \textbf{1.51e+04} \\ \hline
    20 & 76 & \textbf{54} & \textbf{54} & \textbf{3.25e+04} & 1.46e+05 \\ \hline
    20 & 30 & 25 & \textbf{27} & 3.66e+04 & \textbf{1.92e+04} \\ \hline
    28 & 360 & 194 & \textbf{195} & 5.20e+04 & \textbf{2.74e+04} \\ \hline
    28 & 121 & \textbf{84} & \textbf{84} & 7.09e+04 & \textbf{5.96e+04} \\ \hline
    28 & 171 & 108 & \textbf{109} & \textbf{2.74e+04} & \textbf{2.74e+04} \\ \hline
    28 & 300 & \textbf{175} & \textbf{175} & 1.48e+05 & \textbf{1.08e+05} \\ \hline
  \end{tabular}
\end{table}

The Max-Cut results differ from the RSA and synthetic high-order cases. Since the standard Max-Cut objective is naturally quadratic, the original QUBO formulation does not require auxiliary variables and is often competitive. The HUBO version obtained by Algorithm~\ref{alg:Max-Cut_hubo} can still match QUBO quality in many instances, but it does not provide the same structural advantage as in problems whose native model is already high-order. This supports the main interpretation: the preferred formulation depends on whether a transformation reduces or increases the effective dimension faced by the tensor sampler.

\section{Discussion}
\label{sec:discussion}

Our theoretical and empirical investigation reveals a clear and consequential trade-off between HUBO and QUBO formulations within the context of tensor sampling optimization. The central finding is that for problems with inherent high-order interactions, the \textbf{HUBO formulation is not just more expressive but often significantly more efficient.}

\textbf{The Computational Advantage of HUBO:} The prevailing assumption that QUBO is universally more tractable is challenged by our results. The HUBO-to-QUBO reduction process introduces a large number of auxiliary variables, drastically increasing the \textbf{effective dimensionality} of the optimization problem. Since the cost of tensor sampling methods scales sharply with dimensionality (see Appendix~\ref{app:complexity}), this expansion often nullifies any theoretical advantage of the quadratic form. As our experiments show, for problems like RSA factorization and high-order synthetic functions, the native HUBO formulation consistently found better solutions (Tables~\ref{tab:3}, \ref{tab:4}, and \ref{tab:5}) and, in several cases, did so with a \textbf{lower computational budget} (shorter runtimes, fewer evaluations). This demonstrates that preserving the native high-order structure creates a simpler, lower-dimensional search space for tensor methods to explore, leading to superior performance.

\textbf{When QUBO Remains Relevant:} The QUBO formulation retains its value for problems that are naturally quadratic (like the standard Max-Cut problem) or where the number of variables is low. In these cases, as shown in Tables~\ref{tab:maxcut_protes} and \ref{tab:maxcut_ttopt}, both formulations can find solutions of similar quality. The choice may then hinge on the maturity and wider availability of QUBO-specific solvers outside the tensor sampling paradigm.

\textbf{The Critical Role of Dimensionality:} The effectiveness of any formulation is dictated by the capabilities of the underlying optimizer. Methods like PROTES and TTOpt are powerful but have strict practical limits regarding the number of variables they can handle (Table~\ref{tab:tensor_limits}). Therefore, the decision between HUBO and QUBO is, in practice, a decision about which formulation produces a tensor of manageable dimensionality. Our work shows that for high-order problems, the native HUBO formulation is the only viable path to a solution, as the QUBO-reduced form often exceeds the feasible dimension for tensor samplers.

\textbf{Limitations:} The experiments are intended to isolate the effect of formulation and dimension growth rather than to exhaustively benchmark all available solvers. The conclusions are therefore most directly applicable to tensor-sampling optimizers and to reductions that introduce many auxiliary binary variables. Other QUBO reductions, penalty choices, or specialized QUBO solvers may change the practical balance for particular applications. Nevertheless, the dimension-growth mechanism identified here remains relevant whenever the optimizer scales strongly with the number of binary variables.

\section{Conclusion}
\label{sec:conclusion}

This work presented a comparative analysis of Higher-Order (HUBO) and Quadratic (QUBO) Unconstrained Binary Optimization models for probability optimization using tensor sampling techniques. We demonstrated that the traditional pathway of reducing complex problems to QUBO can be counterproductive when using modern tensor methods. The introduction of auxiliary variables inflates the problem dimension, often pushing it beyond the feasible scope of samplers like PROTES and TTOpt and resulting in poorer solutions. Our central conclusion is that for problems with inherent high-order interactions, \textbf{the HUBO formulation is often superior: it provides a more accurate representation of the problem structure and, crucially, can lead to more computationally efficient optimization.} This refines the conventional wisdom that QUBO is always the more tractable choice and provides a dimensionality-aware framework for selecting an optimization formulation. Future work should extend the comparison to additional quadratization schemes, larger benchmark families, and hybrid solvers that can exploit both native high-order structure and mature QUBO optimization tools.

\backmatter

\bmhead{Funding}

The study was supported by the Ministry of Economic Development of the Russian Federation (agreement No. 139-10-2025-034 dd. 19.06.2025, IGK 000000C313925P4D0002).

\appendix
\section{Tensor Train for optimization}
\label{app:tt_optimization}
\paragraph{What is Tensor Sampling?}  
Probabilistic distributions are often high-dimensional, and storing them explicitly is computationally infeasible due to the curse of dimensionality. To address this problem, we can approximate the distribution as a high-dimensional tensor and represent it in the Tensor Train (TT) format. In this way, we avoid explicitly storing the full distribution while still being able to efficiently generate random samples (vectors) from it.

\paragraph{TTOpt.}  
The TTOpt \citet{TTOpt} algorithm is designed to efficiently locate the maximum or minimum elements of a high-dimensional tensor $\mathcal{T} \in \mathbb{R}^{N_1 \times \cdots \times N_d}$ that is defined implicitly through a multivariate function $T(\boldsymbol{x})$. 
Unlike gradient-based methods, TTOpt requires no derivatives and is well-suited for black-box, discrete, or non-smooth optimization problems in high dimensions.

This method is based on the following max-volume observation \citet{maxvol}: 
if $\hat{\boldsymbol{T}}$ is an $R \times R$ submatrix of maximal volume (in selected rows and columns), then the maximal (by modulus) element $\hat{T}_{\text{max}} \in \hat{\boldsymbol{T}}$ bounds the absolute maximal element $T_{\text{max}}$ in the full matrix $\boldsymbol{T}$:
\begin{equation}
\hat{T}_{\text{max}} \cdot R^2 \geq T_{\text{max}}.
\end{equation}
Numerical experiments show that even this bound is pessimistic, that is why TTOpt and other maxvol-based methods aim to find the submatrix with a large volume instead of directly searching for the largest element. 

We begin by considering the first unfolding of the tensor, $\boldsymbol{T}_1 \in \mathbb{R}^{N_1 \times (N_2 \cdots N_d)}$, which reshapes $\mathcal{T}$ into a matrix by grouping all but the first mode. Then, we sample $R_1$ random column indices $I_1^{(C)}$, each corresponding to a multi-index $(n_2, \dots, n_d)$, and  construct a submatrix $\boldsymbol{T}_1^{(C)} \in \mathbb{R}^{N_1 \times R_1}$ by sweeping over all values of $n_1 = 1, \dots, N_1$.

Next, TTOpt applies the \texttt{maxvol} algorithm to $\boldsymbol{T}_1^{(C)}$ to select a subset of $R_1$ rows that form a submatrix $\hat{\boldsymbol{T}}_1 \in \mathbb{R}^{R_1 \times R_1}$ of maximal volume. The corresponding row indices $I_1^{(R)}$ are stored and used in the next step.

The algorithm reinterprets the selected $R_1$ rows as a new implicit matrix $\boldsymbol{J}_2 \in \mathbb{R}^{(R_1 N_2) \times (N_3 \cdots N_d)}$, which corresponds to reshaping the partially selected tensor slice. We then sample $R_2$ random columns $I_2^{(C)}$ in this new matrix and compute the submatrix $\boldsymbol{J}_2^{(C)} \in \mathbb{R}^{(R_1 N_2) \times R_2}$ explicitly.

Again, the algorithm applies \texttt{maxvol} to extract $R_2$ row indices $I_2^{(R)}$, which correspond to combinations of the previously selected indices and new values along the second mode. This process is repeated sequentially across all modes $k = 3, \dots, d$, alternating between random column sampling and maxvol-based row selection.

After reaching the last mode, the algorithm proceeds backward from mode $d$ to mode $1$, now treating previously fixed row indices as variables and selecting new column indices. This alternating forward and backward procedure (referred to as \textit{sweeps}) continues until convergence or until a predefined budget of function evaluations is exhausted.

After $T$ sweeps, the algorithm returns the multi-index $(n_1^*, \dots, n_d^*)$ corresponding to the largest observed value of $T(\mathbf{n})$, along with the value itself.

\textbf{Remark:} The convergence of this method, in the case when it is applied to tensors, has not been proven. We can only be sure that the result will improve monotonically.

\paragraph{Optima-TT}
This method \citet{OptimaTT} finds the maximum or minimum of a tensor given in the TT-format. 

The main idea of this algorithm is to associate each tensor element with the values of the probability density function of a random vector $\boldsymbol{\xi} = \{\xi_1, \ldots, \xi_d\}$, where $\xi_i \in \{1,2,\dots,N_i\}$ and $N_i$ is the size of the corresponding mode ($N_i = 2$ in our case). 

To avoid negative values, we need to square the tensor element-wise: 
$p(n_1, \dots, n_d) = C (\mathcal{T}[n_1, \dots, n_d])^2$ for all $n_k = 1, 2, \dots, N_k$ ($k = 1, 2, \dots, d$), 
where $p$ is the probability distribution of $\boldsymbol{\xi}$ and $C$ is a normalization constant. 

In this formulation, finding the maximum modulo element in the tensor $\mathcal{T}$ is equivalent to finding the most probable value of $\boldsymbol{\xi}$.

The algorithm builds the probability distribution, follows the \citet{PDE}, and sequentially samples the vector $\boldsymbol{\xi}$ from the corresponding marginal distribution at each step. 

\begin{equation}
\xi_i \sim \boldsymbol{p_i}(\xi_i|\xi_1, \dots, \xi_{i-1})
\end{equation}

At each step, we keep the $K$ most probable candidates. If in the $l$-th step of the algorithm we have $\{\tilde{\xi}_1, \dots, \tilde{\xi}_{l-1}\}$, then the marginal distribution function $p_l$ is given by the following expression:

\begin{equation}
p_t(\xi_l \mid \bar{\xi}_1, \dots, \bar{\xi}_{l-1}) = 
\sum_{n_{l+1}=1}^{N_{l+1}} \cdots \sum_{n_d=1}^{N_d} p(\bar{\xi}_1, \dots, \bar{\xi}_{l-1}, \xi_l, n_{l+1}, n_{l+2}, \dots, n_d).
\end{equation}

In this case, efficient computation is ensured by the following theorem.

\textbf{Theorem 2:}

Let $\mathcal{T}$ be a tensor with the TT-representation. Let the TT-cores of this representation be orthogonalized. Then the result of the convolution of this tensor with itself by the last $l$ indices ($0 < l < d$) is given by the following explicit expression:
\begin{equation}
\sum_{n_{d-l+1}=1}^{N_{d-l+1}} \sum_{n_{d-l+2}=1}^{N_{d-l+2}} \cdots \sum_{n_d=1}^{N_d} \mathcal{T}[n_1, \dots, n_d] \mathcal{T}[n_1, \dots, n_d] = \left\| \mathcal{G}_1[1, n_1, :] \mathcal{G}_2[:, n_2, :] \cdots \mathcal{G}_{d-l}[:, n_{d-l}, :] \right\|_2^2.
\end{equation}

\paragraph{PROTES.}  
The PROTES \citet{PROTES} (\emph{Probabilistic Optimization with Tensor Sampling}) method for black-box optimization is based on probabilistic sampling from a probability density function given in the low-parametric TT format. 

This method, following \citet{dolgov2020approximation}, represents the distribution $p(x)$ as a tensor $\mathcal{P}_{\boldsymbol{\theta}} \in \mathbb{R}^{N_1 \times N_2 \times \dots \times N_n}$ (in our case $N_1=\dots=N_n=2$), where high-quality solutions should receive high probability. $\mathcal{P}_{\boldsymbol{\theta}}$ is generated in a low-parametric representation with a set of parameters $\boldsymbol{\theta}$, having the same shape as the target tensor $\mathcal{T}$.

This method starts from the random tensor $\mathcal{P}_0$ and, at each iteration of the algorithm, generates $K$ trial points $x_i$ from the current distribution $\mathcal{P}_\theta$, represented in the TT format.

Among them, $k$ points with the smallest values of the objective function $f(x_i)$ are selected.
These candidates have indices $\mathcal{S} = \{s_1,\dots,s_k\}$ in the set $\mathcal{X}_K$ of sampled candidates.

The parameters $\theta$ are updated using gradient ascent, maximizing the logarithm of the likelihood of the selected candidates:

\begin{equation}
        \widehat{L}_{\boldsymbol{\theta}}(\{\boldsymbol{x}_{s_1}, \boldsymbol{x}_{s_2}, \dots, \boldsymbol{x}_{s_k}\}) = \sum_{i=1}^k \log \left( \mathcal{P}_{\boldsymbol{\theta}}[\boldsymbol{x}_{s_i}] \right).
\end{equation}

\section{Complexity Analysis}
\label{app:complexity}
\paragraph{TTopt} method has the following asymptotic complexity: 
\begin{equation}
\mathcal{O}\left(T \cdot d \cdot \max_{1 \leq k \leq d} (N_k R_k^2)\right),
\end{equation}
where $T$ is the number of sweeps, $d$ is the dimension, $N_k$ is the size of the $k$-th mode (in our case $N_k=2$ for every $k$), and $R_k$ is the TT rank. 

\paragraph{Optima-TT} method has the following asymptotic complexity:
\begin{equation}
    \mathcal{O}\left(d \cdot K \cdot N \cdot R^2\right),
\end{equation}
where $d$ is the tensor dimension, $K$ is the number of selected candidates for each TT-core, and $N$ and $R$ are the typical mode size and TT-rank, respectively.

\paragraph{PROTES} has the following asymptotic complexity:
\begin{equation}
\mathcal{O}\left( M \cdot d \cdot \left( \frac{k}{K} \cdot k_{gd} \cdot \overline{R}^2 + (\overline{N} + \overline{R}) \cdot \overline{R} + \alpha(N) \right) \right),
\end{equation}
where $M$ is the number of requests to the target function, $d$ is the tensor dimension, $k$ is the number of best candidates with minimal objective value, $K$ is the number of sampled candidates, $k_{gd}$ is the number of gradient ascent steps, $\overline{R}$ and $\overline{N}$ are the average TT-rank and mode size, respectively, and $\alpha(N)$ is the complexity of sampling from the generalized Bernoulli distribution with $N$ outcomes.
\bibliography{sample}
\end{document}